\documentclass[10pt,twoside]{amsart}

\usepackage{amsmath}
\usepackage{amssymb}
\usepackage{amsthm}
\usepackage{bm}
\usepackage{bigints}
\usepackage[mathscr]{euscript}
\usepackage{accents}

\usepackage{graphicx}
\usepackage[table]{xcolor}

\usepackage{tabularx}
\usepackage{array}
\usepackage{booktabs}

\usepackage{enumerate}
\usepackage{parskip}
\usepackage[section]{placeins}
\usepackage[a4paper,top=1.4in,bottom=1.4in,
left=1.25in,right=1.25in]{geometry}

\usepackage{url}
\usepackage{epigraph}

\usepackage[
colorlinks=true,
linkcolor=blue,
citecolor=blue,
urlcolor=blue
]{hyperref}
\usepackage{cleveref}

\newtheorem{Th}{Theorem}

\theoremstyle{definition}

\theoremstyle{remark}

\newtheorem*{thank}{Acknowledgment}

\newcommand{\outline}[1]{\textit{#1}}

\newcommand{\lra}{\longrightarrow}
\newcommand{\inter}{\cap}
\newcommand{\union}{\cup}
\newcommand{\un}{\underline}
\newcommand{\til}{\widetilde}
\newcommand{\ov}{\overline}
\newcommand{\wh}{\widehat}
\newcommand{\ot}{\otimes}

\newcommand{\eps}{\varepsilon}
\newcommand{\vep}{\varepsilon}
\newcommand{\va}{\varphi}

\newcommand{\la}{\lambda}
\newcommand{\si}{\sigma}

\newcommand{\egdef}{:=}
\newcommand{\setdef}{\stackrel{\mathrm{def}}{=}}

\newcommand{\ds}{\displaystyle}
\newcommand{\script}{\scriptstyle}
\newcommand{\spa}{\operatorname{span}}

\newcommand{\scalar}[2]{(#1\mid #2)}
\newcommand{\tend}[3][]{\xrightarrow[#2\to#3]{#1}}

\newcommand{\spdj}{\perp_{\mathrm{sp}}}

\DeclareMathOperator{\Id}{Id}
\DeclareMathOperator{\Var}{Var}

\newcommand{\ca}{\mathcal{A}}
\newcommand{\cb}{\mathcal{B}}
\newcommand{\cc}{\mathcal{C}}
\newcommand{\cd}{\mathcal{D}}
\newcommand{\ce}{\mathcal{E}}
\newcommand{\cf}{\mathcal{F}}
\newcommand{\cg}{\mathcal{G}}
\newcommand{\ch}{\mathcal{H}}
\newcommand{\cm}{\mathcal{M}}
\newcommand{\cp}{\mathcal{P}}
\newcommand{\cs}{\mathcal{S}}
\newcommand{\ct}{\mathcal{T}}

\newcommand{\M}{\mathcal{M}}
\newcommand{\BC}{\mathcal{B}}
\newcommand{\Ca}{\mathcal{C}}
\newcommand{\F}{\mathcal{F}}
\newcommand{\J}{\mathcal{J}}
\newcommand{\Lit}{\mathcal{L}}
\newcommand{\Pa}{\mathcal{P}}
\newcommand{\G}{\mathcal{G}}

\newcommand{\xbm}{(X,\mathcal{B},\mu)}
\newcommand{\zdr}{(Z,\mathcal{D},\rho)}
\newcommand{\ycn}{(Y,\mathcal{C},\nu)}

\newcommand{\bs}{\mathbb{S}}
\newcommand{\bd}{\mathbb{D}}
\newcommand{\A}{\mathbb{A}}
\newcommand{\Q}{\mathbb{Q}}
\newcommand{\R}{\mathbb{R}}
\newcommand{\Pro}{\mathbb{P}}
\newcommand{\T}{\mathbb{T}}
\newcommand{\C}{\mathbb{C}}
\newcommand{\Z}{\mathbb{Z}}
\newcommand{\N}{\mathbb{N}}
\newcommand{\E}{\mathbb{E}}
\newcommand{\1}{\mathbb{I}}

\newcommand{\mob}{\boldsymbol{\mu}}
\newcommand{\bml}{\bm{\lambda}}

\newcommand{\Supp}{\Sigma}

\definecolor{headblue}{RGB}{31,56,100}
\definecolor{rowgrey}{RGB}{238,242,248}

\newcolumntype{L}{>{\raggedright\arraybackslash}X}
\newcolumntype{C}[1]{>{\centering\arraybackslash}p{#1}}

\title[O\lowercase{pen}AI partially disproves Littlewood's 1961 $L^\alpha$-flat theorem ]{
 O\lowercase{pen}AI partially disproves Littlewood's 1961 $L^\alpha$-flat theorem}

\author{\lowercase{el} H\lowercase{oucein} \lowercase{el} A\lowercase{bdalaoui}$^\star$}
\address{$^\star$ University of Rouen Normandy ,
	Department of Mathematics, LMRS  UMR 60 85 CNRS\\
	Avenue de l'Universit\'e, BP.12
	76801 Saint Etienne du Rouvray - France .}
\email{elhoucein.elabdalaoui@univ-rouen.fr}

\urladdr{http://www.univ-rouen.fr/LMRS/Persopage/Elabdalaoui/}

\date{\today}
\dedicatory{}

\begin{document}

\maketitle	
\begin{minipage}{0.48\textwidth}
	\epigraph{Everyone writes and nobody reads. \cite{Er}}{P. Erdös (Should be attributed to L. Fejér){\footnote {Does Erdös-Fejér's quotation hold for machines?!}}}
\end{minipage}\hfill
\begin{minipage}{0.48\textwidth}
	\epigraph{A mathematician is a machine for turning coffee into theorems}{P. Erdös. (Probably due to Rényi){\footnote {AI, it seems, is one that no longer needs the coffee.}}}
\end{minipage}\hfill
\begin{minipage}{0.48\textwidth}
	\epigraph{The dangerous machine is not the one that breaks down, but the one that keeps running under corrupt principles.}{ Anoynomus}
\end{minipage}
\hfill
\begin{minipage}{0.48\textwidth}  
	\epigraph{In life, unlike chess, the game continues after checkmate.}{Isaac Asimov}
\end{minipage}

	\begin{abstract} 	We present a counterexample, produced by OpenAI, to Littlewood's 1961 theorem on $L^\alpha$-flat polynomials. The theorem asserts that no $L^\alpha$-flat polynomials exist within the class of trigonometric polynomials of the form$$\sum_{m=1}^{n} a_m \cos(m\theta+\alpha_m), \qquad n \in \mathbb{N}^*,$$satisfying$$k\, n^2 \sum_{m=1}^{n} a_m^2 \;\leq\; \sum_{m=1}^{n} m^2 a_m^2,$$where $k$ is an absolute constant (alternatively expressed on the right-hand side as $K = 1/k$). We show that OpenAI's counterexample is contained within the class of hyper-lacunary sequences. Consequently, it offers only a limited refutation, as such sequences fail to satisfy an implicit assumption of Littlewood's theorem. We also discuss broader questions surrounding this counterexample and provide a brief historical overview.
	\end{abstract}

\section*{Introduction}
\maketitle

In a preprints \cite{OAI-U}, \cite{OAI-C}, \cite{AOI-C2} released on September 23, 2026 by  OpenAI related to the ultraflat polynomials conjecture of  Erdös-Newman\footnote{According to \url{https://erdosproblems.ai/wiki/problems/polynomials/E1150}, only the results of \cite{OAI-C} has been verified by Lean .}, a counterexample to $L^\alpha$ flat polynomials Littlewood Theorem is given. This later Theorem goes back to 1961's paper by Littlewood \cite{L}. Therein the author, generalized an unpublished proof of Erdös theorem which assert that  there is no $L^\alpha$-flat polynomials within the class of trigonometric polynomials of the form 
$$\sum_{j=1}^{n}a_m \cos(m t+\phi_m).$$ 
The  $L^\alpha$ flat polynomials Littlewood (1961's Theorem)  can be stated as follows.
\begin{Th}[1961's Theorem of Littlewood \cite{L}]\label{Lit}
	Let $\ds f_n(t)=\sum_{j=1}^{n}a_m \cos(m t+\phi_m)$ and assume that we have
	$$\sum_{m=1}^{n}a_m^2 \leq \frac{K}{n^2} \sum_{m=1}^{n}m^2a_m^2,$$
	for some absolute constant $K$. Then, for any $\alpha>0$ there exists a constant $A(k,\alpha)$ such that
	$$
	\begin{array}{ll}
		\|f_n\|_\alpha \leq \big(1-A(k,\alpha)\big)\|f_n\|_2, & \hbox{ if $\alpha<2$;} \\
		\|f_n\|_\alpha \geq \big(1+A(k,\alpha)\big)\|f_n\|_2, & \hbox{ if $\alpha>2$.}
	\end{array}
	$$	 
\end{Th}

As a consequence, he obtained Erdős's unpublished result as a corollary. Indeed, when $a_m =\pm 1$, the condition is satisfied with $k=\frac{1}{3}$.

At this point, let us emphasize that this example satisfies the following condition, which is implicit in Littlewood's paper \cite{L}.

$$\frac{\sum_{m=1}^{n}a_m}{|a_n|}=n \to \infty \textrm{~~as~~} n \to \infty.$$
  
The principal ingredients of the proof given by Littelwood are Bernstein-Zygmund inequality and the fact that for any real function $f$
$$|f(\theta)-f(\theta')|=\pm|f(\theta)|\pm |f(\theta')|.$$
Therefore, the relation stated and referred to as ``Littlewood's argument" in \cite{eA2} should be read as.

$$|\Re(P(z))-\Re(P(z'))|=\pm|\Re((P(z))|\pm |\Re(P(z'))|,$$
and 
$$|\Im(P(z))-\Im(P(z'))|=\pm|\Im(P(z))|\pm |\Im(P(z'))|,$$

where $\Re(P)$ is the real part of $P$ and $\Im(P)$ its imaginary part. This fix the issue raised by OpenAI's paper. \\

On the other hand, OpenAI's proof is released as a three papers \cite{OAI-U}, \cite{OAI-C}, \cite{AOI-C2}. But, it is seems to us that the principal paper is the companion paper \cite{OAI-C}. Indeed, therein, OpenAI state that
\begin{Th}\label{Oup}For every $\eta>0$, for large $N$, some polynomial  $P$ with $N$ consecutive plus-minus-one coefficients has modulus between,
	 $\sqrt{N}/16$ and $(1+\eta) \sqrt{N}$ on the whole unit circle, including $z = 1$ and $z = -1$, that is,
	 $$\frac{1}{16} \sqrt{N}\leq \max_{|z|=1}|P(z)| \leq (1+\eta) \sqrt{N}.$$
\end{Th}

From Theorem \ref{Oup} it is easy to see that the sequence of polynomials $(P_N)$ constructed is ultraflat, that is,
\begin{align}
	\max_{\lvert z\rvert=1}\left\lvert \frac{\lvert P_N(z)\rvert}{\sqrt{N}}-1 \right\rvert \to 0
	\qquad \text{as } N\to\infty.
\end{align}
Indeed, since the $L^2$-norm is dominated by the $L^\infty$-norm, we
get that, for large $N$,
\begin{align}
	\left\lVert \frac{P_N(z)}{\sqrt{N}} \right\rVert_2 =1\leq
	\max_{\lvert z\rvert=1}\left\lvert \frac{P_N(z)}{\sqrt{N}} \right\rvert \leq (1+\eta).
\end{align}
Letting $N\to\infty$, we get
\[
\max_{\lvert z\rvert=1}\left\lvert \frac{P_N(z)}{\sqrt{N}} \right\rvert \to 1.
\]
Whence, for any $\eta>0$, for large $N$,
\begin{align}
	(1-\eta)\leq
	\max_{\lvert z\rvert=1}\left\lvert \frac{P_N(z)}{\sqrt{N}} \right\rvert \leq (1+\eta).
\end{align}
The lower bound $\frac{1}{16}$ seems to be obtained by using the Main
Partial Colouring Lemma due to Lovett and Meka \cite[Th.~4]{LM}. This lemma was also used by P.~Balister \& al.,\ in their astonishing paper
\cite{BBMST}, where they answer affirmatively Problem~297 of Erd\H{o}s \cite[p.~297]{E2}. The author gives the explicit constants $c_1=2^{-160}$ and 
$c_2=2^{12}$ \cite[Theorem 2.1]{BBMST}. They further noticed that ``the constants could be improved in Theorem 2.1 somewhat,$\cdots$''.

\section{Some historical facts \protect\footnote{All quotations are given verbatim, as they appear in the cited sources.}}

The ultraflat polynomial problem was initiated by Erd\H{o}s and Newman.
Precisely, P.~Erd\H{o}s wrote about his ultraflat conjecture:
\begin{quote}
	Some of these questions may not be ``serious'' mathematics, but I am
	sure the following final problem, considered by D.~J.~Newman and myself
	for a long time, is both difficult and interesting: Let $\epsilon_k=\pm 1$.
	Is it true that there is an absolute constant $c$ so that, for every
	choice of the $\epsilon_k$'s,
	\[
	\max_{|z|=1}\Big|\sum_{k=1}^{n}\epsilon_k z^k\Big| > (1+c)\sqrt{n}\,?"
	\]
\end{quote}
Latter, in his 1957's paper, Erdös stated the problem as "Problem 22.'' \cite[p.25]{E2}. He also point out that he have 
unpublished proof that if the polynomial are of the form
$$P(\theta)=\sum_{m=1}^{n}\varepsilon_n\cos(n\theta), \varepsilon=\pm 1, m=1,\cdots, n.$$
Then 
$$\max_{0<\theta<2\pi }\big|P(\theta) |>(1+c)\sqrt{n/2}.$$

The coefficients $(a_k)$ of the polynomials in Problem 22 are unimodular ($|a_k|=1$). Subsequently, in 1961's, Littlewood stated and give a proof of his 1961's $L^\alpha$-flat Theorem \cite{L}.  Let us also point out that P. Erdös in \cite[p.297-298]{E2} asked on whether there exists a sequence of analytic polynomial
$P_n$ with coefficients $\epsilon_{k,n} \in \{\pm 1\}$, degree $n$ and a positive constant $c_1,c_2$ such that
$$c_1 \sqrt{n}< \max_{z: |z|=1}\big|P_n(z)\big| < c_2 \sqrt{n}.$$
He mentioned also that Cluny has a partial unpublished result.Nowadays, such sequences of polynomials are called ultraflat polynomials. 
We summarizes in Table \ref{Tab1} the other contributions from 1957 to 1995.
While a comprehensive history of this problem would be of great interest, the brief nature of this paper requires that it be deferred to a forthcoming monograph 
on flat polynomials.
We therefore restrict ourselves to mentioning that Andrew Odlyzko wrote in \cite{O}:
	``Extensive computations described in this paper\footnote{This referred to \cite{O}.} strongly suggest such polynomials do not exist,
''
\begin{table}[ht]
	\centering
	\small
	\renewcommand{\arraystretch}{1.5}
	\rowcolors{2}{rowgrey}{white}
	\begin{tabularx}{\linewidth}{@{} >{\bfseries}C{3.2cm} C{1.3cm} L @{}}
		\rowcolor{headblue}
		\toprule
		\textcolor{white}{\textbf{Authors}} & \textcolor{white}{\textbf{Year}} & \textcolor{white}{\textbf{Contributions}} \\
		\midrule
		Beller--Newman & 1974 &
		There exist polynomials $P_n$ with $|\widehat{P_n}(k)|<1$ such that
		$\min_{|z|=1}|P_n(z)| \geq C\sqrt{n}$. \\
		
		Byrnes & 1977 &
		For every $\delta>0$ there exist polynomials $P_n$ with $|\widehat{P_n}(k)|=1$
		such that $\operatorname{mes}\{\theta : (P_n) \text{ is flat in the Littlewood sense}\}>1-\delta$. \\
		
		K\"{o}rner & 1979 &
		By introducing random arguments into Byrnes' construction, he showed that
		Littlewood's question has a positive answer. \\
		
		Kahane & 1979 &
		Using a random lemma of K\"{o}rner, he produced an ultraflat sequence $(P_n)$
		whose coefficients $\widehat{P_n}(k)$ are complex numbers of modulus $1$ and
		$(1-\epsilon_n)\sqrt{n} \leq |P_n(z)| \leq (1+\epsilon_n)\sqrt{n}$,
		where $\epsilon_n\to 0$. \\
		
		Saffari--Smith & 1988 &
		They claimed to have proved the $L^\infty$ Erd\H{o}s conjecture, but two years
		later they published a note retracting the claim. \\
		
		Beck & 1990 &
		Using Kahane's construction, he produced flat polynomials in the Littlewood
		sense with coefficients that are roots of unity of order $400$. Using Byrnes'
		construction, he also showed that the coefficients can be taken to be roots of
		unity of order $3$. \\
		
		Benedetto and Hui-Chuan & 1995 &
		They showed that Byrnes' proof contains a flaw that cannot be repaired.
		This does not affect Kahane's result, since Kahane uses Gaussian polynomials (that is, with coefficients $\big(e^{i \pi \frac{k^2}{n}}\big)$
		rather than Byrnes' polynomials. Historically, Gaussian coefficients were first used by D. J. Newman in 1965 to construct $L^1$-flat polynomials \cite{N65} \\
		\bottomrule
	\end{tabularx}
	\caption{Historical milestones on flat polynomials.}\label{Tab1}
\end{table}

\section{OpenAI's counter example.} In \cite{OAI-C}, the coefficients $(a_m)$ are choosing as a hyper-lacunary sequence. We remind that a positive sequence  $(a_m)$ is said to be lacunar if for $m \in \N$, we have
$$\frac{a_{m+1}}{a_m} \geq q >1.$$ 
and hyper-lacunar when
$$\frac{a_{m+1}}{a_m} \to +\infty ~~\textrm{~~as~~} m \to \infty.$$ 
In the OpenAI's example, the sequence is hyper-lacunary, given by $a_m=2^{2^m}$. We thus obtain that the condition in Theorem \ref{Lit} is satisfied with
$k=1/2$. Furthermore, asymptotically, one coefficient dominated all the coefficients. Indeed,
$$\frac{\sum_{m=1}^{n}a_m}{a_n}\tend{n}{+\infty} 1.$$
This give, that for any $z, |z|=1$,
$$\Big|\frac{\sum_{m=1}^{n}a_m z^m}{\|P_n\|_2}-\frac{a_n}{\|P_n\|_2}z^n\Big| \tend{n}{+\infty} 0.$$
We deduce that the sequence is ultraflat, since obviously the trivial
sequence $(\pm 1.z^n)$ is ultraflat.

This strategy is well-known in the Riesz product theory \cite{Z}, \cite{Abd-nad3}, \cite{Abd-Nad-etds}. 

The example does not take the signs of the coefficients into account. However, the strategy remains valid for any sequence $(a_n)$ whose sequence of absolute values $(|a_n|)$ is lacunary, regardless of the signs of the coefficients.

Indeed, consider the polynomials $P_N$ normalized in the $L^2$-norm, that is,
$$P_N(z)=\frac{\sum_{m=1}^{N}a_mz^m}{\|P_N\|_2}.$$
It follows that $(P_N(z))$ are ultraflat. Indeed, we have
\begin{align}
	\left\lVert \frac{\sum_{m=1}^{N}a_mz^m}{|a_N|} \right\rVert_2 \leq
	\max_{\lvert z\rvert=1}\left\lvert \frac{\sum_{m=1}^{N}a_mz^m}{|a_N|} \right\rvert \leq 
	\frac{\sum_{m=1}^{N}|a_m|}{|a_N|}.
\end{align}
Letting $N\to\infty$,  we get 

$$	\max_{\lvert z\rvert=1}\left\lvert \frac{\sum_{m=1}^{N}a_mz^m}{|a_N|} \right\rvert \to 1.$$ 

Another strategy to construct an ultra-flat polynomials was given in \cite{Abd-Nad-etds}. Therein, the authors used the partial rest of the entire series associated to 
the function 
$$\frac{z-\alpha}{1-\alpha .z}, ~~~~~\alpha \in (0,1).$$
This later function is related to the inner-outer functions in $H^p$ Theory. \\

To the best of the author's knowledge, there is a third strategy for producing $L^1$-flat idempotent polynomials (polynomials with coefficients in $\{0,1\}$), due to the author \cite{eA-S}. This later strategy is based on Singer's theorem.

In the next section, we discuss the OpenAI's example and raise a question as to whether it is possible to produce an example that would completely
invalidate 1961's Littlewood's theorem. Although the author expresses some reservations regarding this.

\section*{Discussion.} The author of this note use 1961's theorem of Littlewood \cite{L} combined with \linebreak Clarkson's second inequality to generalize Theorem \ref{Lit}. He  further point out that the condition in Theorem \ref{Lit} is related to Bernstein-Zygmund inequality. Precisely, Littlewood's condition insure that the $L^2$-norm in this class are equivalent in the following sense:
$$k . n .\big\|P_n\big \|_2 \leq \big\|P_n' \big\|_2 \leq  n \big\|P_n\big \|_2.$$
It is well known that the equality holds for $P(z)=c z^n$. The counterexample in \cite{OAI-C} concern the asymptotic case. It does not affect the theorem, if we further require, 
$$\frac{\sum_{m=1}^{n}a_m}{|a_n|} =\frac{\|(a_m)_{m=1}^{n}\|_1}{\|(a_m)_{m=1}^{n}\|_{\infty}} \tend{n}{+\infty}+\infty.$$
This condition takes the signs into account, and it is satisfied by the trivial case $a_m=\varepsilon_m=\pm 1$, $m=1,\cdots,n$. Indeed, 
\[
\frac{\|(a_m)_{m=1}^{n}\|_1}{\|(a_m)_{m=1}^{n}\|_{\infty}}=n \to \infty ~~{as}~~n \to \infty.
\]


At this point, the author would like to complement this short note with the following statement concerning the OpenAI preprints \cite{OAI-U, OAI-C, AOI-C2} on ultraflat polynomials[\footnote{cf. \url{https://github.com/openai/math/blob/main/README.md}}]:

\begin{quote}
	The repository states that its manuscripts were ``produced by an internal OpenAI model'' and that the collection ``includes results at different stages of verification''. Not all manuscripts come with Lean formalizations, and, in the repository's own words, ``Some of the formalized results could have issues''. The README in the release folder of the manuscript gives only the title, the author ``OpenAI'', the date September 23, 2026, and the citation block above; it contains no statement about human assistance. The manuscript names no author other than ``OpenAI'', carries neither an arXiv identifier nor a journal reference, and cites a companion release manuscript (a three-torus diffeomorphism with simple Lebesgue spectrum) for a different realization of its spectral consequence.
\end{quote}

For the Lean formalization, we refer to \url{https://sites.pitt.edu/~hajlasz/OpenAI_Math_Catalogue.html}. There, the manuscripts \cite{OAI-U} and \cite{AOI-C2} are listed as having no Lean formalization. whereas \cite{OAI-C} is listed as having a Lean formalization \footnote{ See also \url{https://erdosproblems.ai/wiki/library/polynomials} and \url{https://www.erdosproblems.com/1150}, which report conflicting statuses.}.

 In light of ongoing developments and pending peer-reviewed validation of concurrent preprints in this domain,  the present manuscript and \cite{eA2} will be adapted or withdrawn to reflect the consolidated literature. The author wishes to point out that his paper \cite{eA2} was under revision when OpenAI released three papers on ultraflat polynomials. It is likely that the revised version will be posted on arXiv before the outcome of the review process is known. This may be partly due to the fact that the mathematical community, including researchers using machine assistance, could require some time to assess these three papers and to clarify the status of OpenAI's results on ultraflat polynomials. 
  
 \begin{thank} The author would like to thank Terrence Adams for sending him th1e link to OpenAI's release of October 7, 2026 for bringing to his attention preprints on dynamical systems, in particular the one on "Simple Lebesgue spectrum", and for the subsequent email exchanges on these and related topics. He also thanks Jean-Paul Thouvenot for bringing to his attention the preprints on "ultraflat polynomials", and for discussions on the consequences of this release for researchers in mathematics.

The presentation, the correction of misprints, and some of the arguments have
benefited in part from suggestions generated in discussions with LLM tools
(Claude by Anthropic, Gemini, ChatGPT, and Mistral). The author takes full
responsibility for the mathematical statements and proofs.

 \end{thank}

\end{document}